\documentclass[11pt]{article}
\usepackage[margin=4.46cm]{geometry}
\usepackage{graphicx}
\usepackage{epsfig}
\usepackage{amsmath}
\usepackage{amssymb}
\newcommand{\commentout}[1]{}

\newcommand{\ra}{\rightarrow}

\newcommand{\nit}{\noindent}
\newcommand{\no}{\nonumber}
\newcommand{\be}{\begin{equation}}
\newcommand{\ee}{\end{equation}}
\newcommand{\ba}{\begin{eqnarray}}
\newcommand{\ea}{\end{eqnarray}}
\newcommand{\bi}{\begin{itemize}}
\newcommand{\ei}{\end{itemize}}
\newcommand{\br}{\begin{eqnarray}}
\newcommand{\er}{\end{eqnarray}}

\newcommand{\fee}{\mbox{$\varphi$}}

\newcommand{\lam}{\mbox{$\lambda$}}
\newcommand{\eps}{\mbox{$\epsilon$}}

\newcommand{\qed}{\mbox{$\square$}\newline}
\newcommand{\dpart}[2]{\frac{\partial #1}{\partial #2}}

\newtheorem{theo}{Theorem}[section]

\newtheorem{lem}{Lemma}[section]

\newtheorem{rmk}{Remark}[section]

\begin{document}

\title{{\Large\bf{Asymptotics for turbulent flame speeds\\ of the viscous G-equation \\
enhanced by cellular and
shear flows}}}
\author{Yu-Yu Liu$^{1}$,\;
Jack Xin$^{2}$, \; Yifeng Yu$^{3}$ \thanks{$^{1,2,3}$Department of Mathematics,
UC Irvine, Irvine, CA 92697, USA.}}

\date{}
\maketitle
\begin{abstract}
G-equations are well-known front propagation models in turbulent combustion and describe the front
motion law in the form of local normal velocity equal to a constant (laminar speed) plus the
normal projection of fluid velocity. In level set formulation, G-equations are
Hamilton-Jacobi equations with convex ($L^1$ type) but non-coercive Hamiltonians.
Viscous G-equations arise from either numerical approximations or regularizations by small
diffusion. The nonlinear eigenvalue $\bar H$ from the cell problem of
the viscous G-equation can be viewed as an approximation of the inviscid
turbulent flame speed $s_T$.  An important problem in turbulent combustion theory
is to study properties of $s_T$, in particular how $s_T$ depends on the flow
amplitude $A$. In this paper, we will study the behavior of $\bar H=\bar H(A,d)$
as $A\to +\infty$ at any fixed diffusion constant $d > 0$.  For the cellular flow,
we show that
$$
\bar H(A,d)\leq O(\sqrt {\mathrm {log}A})  \quad \text{for all $d>0$}.
$$
Compared with the inviscid G-equation ($d=0$),
the diffusion dramatically slows down the front propagation.
For the shear flow,  the limit

\nit $\lim_{A\to +\infty}{\bar H(A,d)\over A} = \lambda (d) >0$ where
$\lambda (d)$ is strictly decreasing in $d$, and has zero derivative at $d=0$. The linear growth
law is also valid for $s_T$ of the curvature dependent G-equation in shear flows.
\end{abstract}

\hspace{.1 in} {\bf Key Words}: viscous G-equations, cellular flows, speed

\hspace{.1 in} enhancement, absence of power law, shear flows, linear law.

\medskip

\hspace{.1 in} {\bf AMS Subject Classification:} 70H20, 76M50, 76M45, 76N20.
%\bigskip

%\hspace{.1 in} {\bf Short Title}: Front Speed Asymptotics of viscous G-equations.

\medskip
\newpage

\section{Introduction}
\begin{figure}[ht] \label{mod}
\centering
\centerline{\includegraphics[width=4in,height=2.0in]{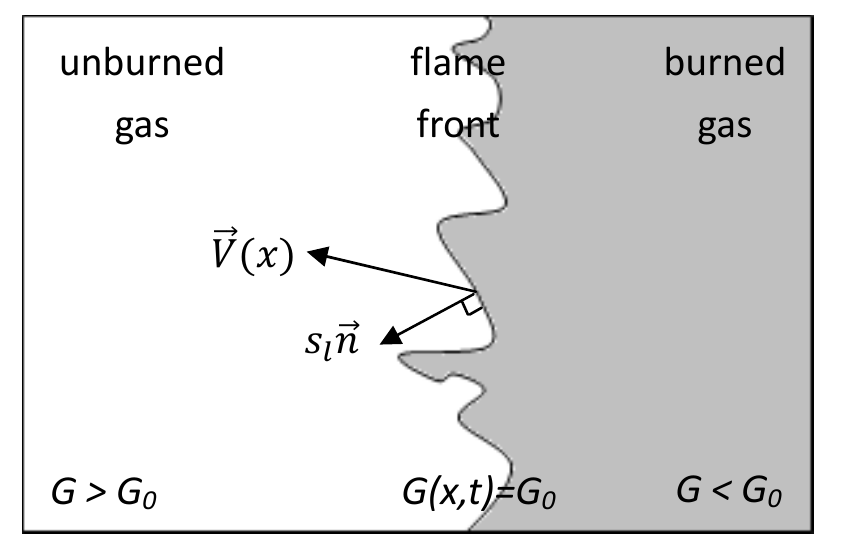}}
\caption{Illustration of local interface velocities in the G-equation and a flame front.}
\end{figure}

\setcounter{equation}{0}
The G-equation has been a very popular field model in combustion and physics literature
for studying premixed turbulent flame propagation \cite{Mark,W85,Siv,S89,Yak,Pet92,EMS,Chert_98,Pet00,Ab_02,OF,Xin_09}.
The inviscid G-equation
on a flame moving in a steady flow has the following form:
\be
G_t+V(x)\cdot DG+s_l|DG|=0, \label{Ge1}
\ee
where $G$ is the level set function of the flame, $V$ is the ambient
fluid velocity field, and the positive constant $s_l$ is called {\it laminar flame speed}.
The constant $s_l$ describes how fast the flame propagates when the fluid is at rest.
The G-equation can be derived through the level set method based on a simple motion law:
the flame propagation speed along the normal direction is equivalent to
$s_l$ plus the normal projection of the fluid velocity (see Figure \ref{mod}).
The level set $\{(x,t): G(x,t)=G_0\}$ of the solution $G=G(x,t)$ represents the
flame front at time $t$. We assume that the flow field $V$ is periodic and incompressible.

Suppose that the initial flame front is planar and
the flame is propagating in the direction $P$ ($G(x,0)=P\cdot x$).
Due to the movement of the fluid,
the flame front will be wrinkled in time.
Eventually, the front will evolve into an asymptotic state moving at a constant speed $s_T$
which depends on $P$ and is called ``turbulent flame speed" in combustion literature.
It can be computed as $s_T=-{\lim_{T\to +\infty}{G(x,T)\over T}}$.
The $s_T$ is conjectured to exist even when $V$ is stochastic and is used in the
combustion community to describe the average speed of a fluctuating front \cite{Pet92}.
To predict and analyze properties of $s_T$ is a fundamental problem in
turbulent combustion theory. When $V$ is periodic in space, the $s_T$ can be
studied in the framework of the periodic homogenization theory of
Hamilton-Jacobi equation \cite{LPV,Evans92}.
It is same as the effective Hamiltonian $\bar H(P)$ of a nonlinear eigenvalue problem
(so called cell problem):
\be
s_l|P+Dw|+V(x)\cdot (P+Dw)=\bar H(P). \label{Ge2}
\ee
Due to the lack of coercivity of the Hamiltonian of the G-equation, the periodic homogenization
and the existence of $\bar H(P)$ have been rigorously
established only very recently by two of the authors \cite{XY_10}
and Cardaliaguet-Nolen-Souganidis \cite{CNS} independently.
When $n=2$, Nolen and Novikov \cite{NN} proved the existence of $\bar H(P)$
for stationary ergodic flows.

In computation of the hyperbolic equation (\ref{Ge1}), certain amount of
numerical diffusion is often present as in Lax-Friedrichs type schemes \cite{OF}.
On the other hand, it is known that as $t$ gets large,
the level set $\{G(x,t)=0\}$ might become quite irregular and cause numerical difficulties.
Among the various regularizations to fix this problem,
one way is to add a diffusion term \cite{FIL} to (\ref{Ge1}) which leads to the
following viscous G-equation
\be
-d\Delta G+G_t+V(x)\cdot DG+s_l|DG|=0, \label{Ge3}
\ee
${d\over s_l}>0$ is called the Markstein diffusivity.
If we consider $-{\lim_{T\to +\infty}{G(x,T)\over T}}$,
the limit $\bar H(P,d)$ is given by  the cell problem
\be
-d\Delta w+s_l|P+Dw|+V(x)\cdot (P+Dw)=\bar H(P,d). \label{Ge4}
\ee
The existence of $\bar H(P,d)$ and classical solutions to (\ref{Ge4}) (unique up to a constant)
can be easily deduced from the standard elliptic regularity theory.
$\bar H(P,d)$ can be viewed as an approximation of the turbulent flame speed $s_T$.

A central issue we address here is the comparison of qualitative behavior of
$\bar{H}(P,d)$ and $\bar{H}(P)$ as we vary $d$ and the amplitude of the flow field.
%Besides the existence, an important question about $\bar H$ is
%its dependence on the fluid velocity.
To this end, let us scale $V$ to $A\, V$ for some positive constant $A$ (flow intensity), so:
\be
-d\Delta w+s_l|P+Dw|+AV(x)\cdot (P+Dw)=\bar H(P,A,d). \label{Ge5}
\ee
An interesting question is to figure out how $\bar H$ behaves as
a function of $A$. There are a few results in the combustion literature
in this direction when $A$ is small, see the references of \cite{R}.
In this paper, we are interested in the asymptotic behavior of $\bar H(P,A,d)$ as $A\to +\infty$.
The usual inf-max formula
\be\label{inf}
\bar H(P,A,d)=\inf_{h\in C^2(\Bbb T^n)}\max_{\Bbb T^n}(-d\Delta h+s_l|P+Dh|+AV(x)\cdot (P+Dh))
\ee
only provides that $\bar H(P,A,d)\leq C(d) \,(A+1)$ which is in general too rough.
Experimental studies show that the turbulent flame speed may grow slower than linear
in some situations. This is the so called ``bending effect", see \cite{S89,R} among others.

The paper is organized as follows. In section 2, we look at the case when $V$ is a two dimensional
cellular flow, $V=\nabla^{\perp} {\cal H} \equiv \nabla^{\perp} \sin 2\pi x_1 \sin 2\pi x_2$.
It is known that $\bar H(P,A,0)=O({A\over \log A})$ for the inviscid case.
The ``bending effect'' occurs marginally.  When the diffusion is large ($d\gg 1$),
it is proved in \cite{NXY} that $\bar H(P,A,d)$ drops dramatically and
has an upper bound as $\sqrt {\log A}$. In the small diffusion regime ($d\ll 1$),
the analysis becomes much more subtle
since the nonlinear $L^1$ term begins to compete with the linear diffusion term.
By a novel ${\cal H}$-weighted gradient estimate of solutions of the cell problem,
we establish the $\sqrt {\log A}$ upper bound for any positive diffusivity $d > 0$.
 Precisely speaking,
the main result of this section is that for a positive $d$-dependent constant $C(d)$
\be
\bar H(P,A,d)\leq C(d) \sqrt {\mathrm {log}A},  \quad \text{for all $d>0$}, \; \text{$A\geq 2$}.
\label{Ge6}
\ee
Our estimates reveal the retention of positive mass of the gradient of solution
in the boundary layers as $A \ra +\infty$, or the loss of the gradient mass in the interior of
each of the quarter cells (Fig. \ref{bd_p1}). The constant $C(d)$ depends on the percentage of such a mass loss.
Our analysis suggests a similar lower bound for $\bar H$ in (\ref{Ge6}) if the $L^1$ norm of the gradient of
a linear advection-diffusion equation arising in advection enhanced diffusion problem \cite{NPR}
satisfies a square root logarithm lower bound which we conjecture to hold.
In essence, we have overcome the obstacle of nonlinearity. As $d\to 0$, the constant $C(d)$ will blow up. 

In section 3, we study the case when $V$ is a shear flow.
We prove that the limit $\lim_{A\to +\infty} {\bar H(P,A,d)\over A}$ is a positive
constant (no ``bending effect")  and is strictly decreasing
with respect to the diffusivity $d$.
The converging rate as $d\to 0$ is also discussed.
Our approach can be used to recover an earlier result in \cite{JKM} about
the forced Burgers' equation. We also investigate the limit for
the curvature dependent G-equation, i.e, replace the diffusion term
by the mean curvature of the flame front. We showed that the limit is the same
as for the inviscid G-equation. Our results in this section
are consistent with the natural intuition on the front propagation speed:
\be \label{speed_comp}
\mathrm{Viscous\ speed}\leq \mathrm{Curvature\ dependent\ speed }\leq \mathrm {Inviscid\ speed}.
\ee

We remark that if the flow field is compressible, the situation is very different.
Firstly, positive diffusivity may increase the propagation speed.
Secondly, there may be flame trapping (in the inviscid case) or
exponential decay of front speed (in the viscous case) due to
the high turbulence intensity ($A \gg 1$).
Explicit analytical results of this sort for the one space dimensional G-equations are
reported in \cite{LXY}.

In section 4, we show numerical results of $\bar{H}$ in viscous
G-equations and cellular flows, and propose an empirical law
$\bar{H}(P,d,A) \sim c(d) \sqrt{\log A}$, $A\gg 1$, $d> 0$ fixed, $c(d)$ decreasing in $d$.
The paper ends with concluding remarks in section 5.
Section 6, the appendix, provides a proof of $L^p$ ($p\in [1,2]$) gradient estimate of
a linear advection-diffusion equation, which is needed in proving the main results.

The work was partially supported by NSF grants
DMS-0712881 (JX) and DMS-0901460 (YY).

\section{Root Log Upper Bound in Cellular Flows}
\setcounter{equation}{0}
Without loss of generality, we assume that $s_l=1$.  Consider the cell problem
\be\label{bd1}
d\Delta w+|P+Dw| + AV(x)\cdot (P+Dw) = \bar H(P,A,d).
\ee
Here we switch $-d\Delta w$ to $d\Delta w$ through a simple change of variables.
In this section, let us look at front speeds in cellular flows.  A typical example is
$V =\nabla^{\bot}\cal{H}$ where the stream function $\cal{H}$$(x)=\sin(2\pi x_1)\sin(2\pi x_2)$.
For simplicity, we will work with this example and write $\cal{H}$ as $H$ hereafter.
The following is the main result of this section which says that the
(viscous) turbulent flame speed increases no faster than square root of $\log A$.

\begin{figure}\centering\label{bd_p1}
\includegraphics{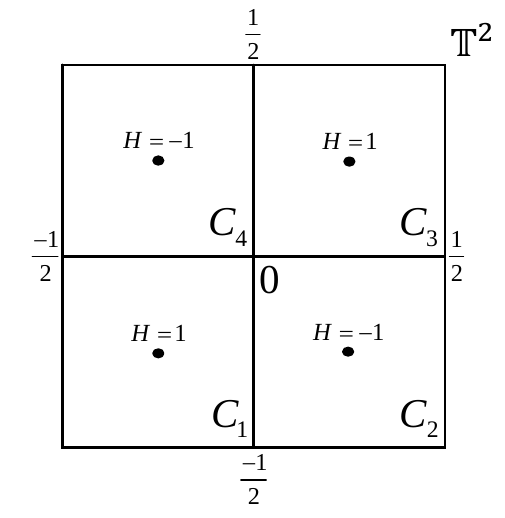}
\caption{Decomposition of a unit cell into quarter cells $C_i$, $i=1,2,3,4$.}
\end{figure}

\begin{theo}\label{bd_t1}
\[
\bar H(P,A,d)\leq O(\sqrt {\mathrm {log}A})  \quad \text{for all $d>0$}.
\]
\end{theo}
Let $e_1=(1,0)$ and $e_2=(0,1)$.  Proofs for $P=e_1$ and $P=e_2$ are similar. Also, the inf-max formula (\ref{inf}) implies
that $\bar H$ is a convex and homogeneous of degree one as function of $P$.
Hence it suffices to prove the above theorem for $P=(1,0)=e_1$.
Let us denote $\bar H(e_1,A,d)=\lam_A$, and omit $d$ dependence for the moment.
Clearly $1\leq \lam_A \leq C(A+1)$. Hereafter, $C$ denotes a constant
independent of the flow intensity $A$.  Note that $C$ might depend on the diffusivity constant $d$.  We also assume that $A\geq 2$.
Also we split $\Bbb T^2=[-{1\over 2},{1\over 2}]\times [-{1\over 2},{1\over 2}]$ into
four cells $C_1$, $C_2$, $C_3$, $C_4$ shown as in Figure \ref{bd_p1}.

Let $G=e_1\cdot x+w(x)$ and $\int_{\Bbb T^2}Gdx=0$. Then
\be\label{bd2}
d\Delta G + |DG|+AV(x)\cdot DG=\lam_A=\bar H(e_1,A,d).
\ee
Integrating both sides and using the incompressibility of $V$ and $V\cdot DH=0$, we have that
$$
\int_{\Bbb T^2}|DG|\,dx=\lam_A.
$$
Let us denote $v_A={G\over \lambda_{A}}$. Then
\be\label{2e.2}
d\Delta v_A+|Dv_A|+AV(x)\cdot Dv_A=1
\ee
and
\be\label{L1bound}
\int_{\Bbb T^2}|Dv_A|\,dx=1.
\ee
Since $\int_{\Bbb T^2}v_A\,dx=0$, we have that
\be\label{mass}
||v_A||_{W^{1,1}(\Bbb T^2)}\leq C.
\ee
Owing to the Sobolev inequality,
$$
\int_{\Bbb T^2}v_{A}^{2}\,dx\leq C.
$$
Upon a subsequence, we may assume that
\[
v_A \rightharpoonup v\ \ \mbox{in}\ \ L^2(\Bbb T^2)\]

We first proves several lemmas. The first one says that the $H^1$ norm of $v_A$ is locally bounded.
\begin{lem}\label{2.2}
$$
\int_{\Bbb T^2}|Dv_A|^2 H^2\,dx\leq C.
$$
\end{lem}

\nit {\bf Proof}: Multiply   (\ref{2e.2}) by $H^2v_A$.  The incompressibility of $V$ and $V\cdot DH=0$
imply that
$$
\int_{\Bbb T^2}V\cdot Dv_Av_AH^2\,dx=0.
$$
Then using integration by parts and Cauchy's inequality, we get that
$$
\int_{\Bbb T^2}|Dv_A|^2H^2\,dx\leq C(1+\int_{\Bbb T^2}v_{A}^{2}\,dx)\leq C.
$$
So the above lemma holds.
\qed

\begin{rmk}\label{boundofv}
From the above lemma, it is clear that for all $\epsilon>0$,  $v\in H^1(\{|H|>\epsilon\})$ and
$$
\int_{\{|H|>\epsilon\}}|Dv|^2 H^2\,dx\leq C.
$$
Then Sobolev embedding implies that
$$
\lim_{A\to +\infty}||v_A-v||_{L^2(\{|H|>\epsilon\})}=0.
$$
\end{rmk}

Next we show that the oscillation of  $v_{A}$  along nonzero level curves of $H$ will tend to zero.
\begin{lem}\label{levelcurve}
$$
\int_{\Bbb T^2}H^4|V(x)\cdot Dv_A|^2\,dx\leq {C\over A}.
$$
In particular,
\be\label{curve}
V(x)\cdot Dv=0  \quad \text{for a.e  $x\in \Bbb T^2$}.
\ee

\end{lem}

\nit {\bf Proof}:  To ease notation, we write $v_A=U$ in this proof.  Then
$$
d\Delta U+|DU|+AV(x)\cdot DU=1.
$$
Multiplying  $H^4V(x)\cdot DU$ on both sides of the above equation and integrating over
$\Bbb T^2$ show that
\[
d\int \Delta U(V(x)\cdot DU)H^4dx + \int |DU|(V(x)\cdot DU)H^4dx +A\int (V(x)\cdot DU)^2H^4dx
\]\[
=\lam_A \int V(x)\cdot DUH^4\, dx =0
\]
The last equality is due to the incompressibility of $V$ and $V\cdot DH=0$.  Note that
\[\int \Delta U(V(x)\cdot DU)H^4dx =
\int \left(\sum_{i=1}^2U_{x_ix_i}\right)\left(\sum_{k=1}^2V_kU_{x_k}\right)H^4dx
\]
where $V=(V_1,V_2)$ and integration by parts,
\[
=-\underbrace{\int \sum_{i=1}^2\sum_{k=1}^2H^4U_{x_i}V_{kx_i}U_{x_k}\,dx}_{\textbf{I}}
-\underbrace{\int \sum_{i=1}^2\sum_{k=1}^2H^4U_{x_i}V_kU_{x_kx_i}\,dx}_{\textbf{I\!I}}
\]\[
-\underbrace{4\int_{\Bbb T^2}H^3(DU\cdot DH)(V(x)\cdot DU)\,dx}_{\textbf{I\!I\!I}}
\]
Note that
$$
\textbf{I}+\textbf{I\!I\!I}\leq C\int_{\Bbb T^2}H^2|DU|^2\,dx.
$$
Moreover,
$$
\textbf{I\!I}={1\over 2}\int_{\Bbb T^2} \sum_{k=1}^2V_k\left(|D U|^2\right)_{x_k}H^4 \, dx\, =0
$$
The last inequality is due to the incompressibility of $V$ and $V\cdot DH=0$. Furthermore, Cauchy inequality implies that
$$
|\int |DU|(V(x)\cdot DU)H^4dx|\leq {1\over 2}\left (\int_{\Bbb T^2}|V(x)\cdot DU|^2H^4\,dx+
\int_{\Bbb T^2}|DU|^2H^4\,dx \right ).
$$
Hence Lemma \ref{levelcurve} follows from Lemma \ref{2.2}.
\qed

The next lemma says that $v_A$ converges to $v$ locally in $H^1$ norm.

\begin{lem}\label{bd_l3} Let $W=g(H^2)$ for some $g\in C_{c}^{\infty}((0,1])$. Then
\[
\lim_{A\ra +\infty}\int_{\Bbb T^2}|Dv_A-Dv|^2W\, dx=0
\]
\end{lem}

\nit {\bf Proof}: In fact
\[
\int_{\Bbb T^2}|Dv_A-Dv|^2W\, dx =\int_{\Bbb T^2}(Dv_A-Dv)\cdot(Dv_A-Dv)W\, dx
\]\[
=\underbrace{\int_{\Bbb T^2}Dv_A\cdot(Dv_A-Dv)W\, dx}_{\textbf{I}}
-\underbrace{\int_{\Bbb T^2}Dv\cdot(Dv_A-Dv)W\, dx}_{\textbf{I\!I}}
\]
Due to Lemma \ref{2.2},  $Dv\in L_{loc}^2(\Bbb T^2\backslash \{H=0\})$ and
$Dv_A \rightharpoonup Dv\ \ \mbox{in}\ \ L_{loc}^{2}(\Bbb T^2\backslash \{H=0\})$.
Hence $\mbox{I\!I}\ra 0$ as $A\ra +\infty$. Also,
\[
I = -\int_{\Bbb T^2}\Delta v_A(v_A-v)W\, dx
-\underbrace{\int_{\Bbb T^2} Dv_A(v_A-v)\cdot DW\, dx}_{\textbf{I\!I\!I}}
\]\[
={1\over d}\underbrace{\int_{\Bbb T^2}|Dv_A|(v_A-v)W\, dx}_{\textbf{I\!V}}
+{A\over d}\underbrace{\int_{\Bbb T^2}V(x)\cdot Dv_A(v_A-v)W\, dx}_{\textbf{V}}
-{\textbf {I\!I\!I}}+\textbf{VI\!},
\]
where
$$
\textbf{VI\!}=-{1\over d}\int_{\Bbb T^2}(v_A-v)W\,dx.
$$
By Remark \ref{boundofv} and Lemma \ref{2.2}, $\mbox{I\!I\!I},\mbox{\ I\!V}, \mbox{VI\!} \ra 0$ as $A\ra +\infty$.

Since $V(x)\cdot Dv=0$ and $V\cdot DH=0$,
integration by parts and the incompressibility of $V$ imply that the 5th term {\textbf{V}=0}.
\qed

\begin{rmk}\label{L1convergence}
It follows immediately from Lemma \ref{bd_l3} that for all $\epsilon>0$
$$
\lim_{A\to +\infty}\int_{\{|H|\geq \epsilon\}}|Dv_A-Dv|\,dx=0.
$$
\end{rmk}

\begin{lem}\label{bd_l4}
(i) For $i=1,3$, $v=f_i(H)$ in $C_i$, where $f_i\in C^{2}((0,1))$ and $f_i'< 0$.
(ii) For $i=2,4$, $v=f_i(H)$ in $C_i$, where $f_i\in C^{2}((-1,0))$ and $f_i'> 0$.
\end{lem}

\nit {\bf Proof:} Since $v\in H_{loc}^{1}(\Bbb T^2\backslash\{H=0\})$ and $V(x)\cdot Dv=0$ a.e.,
it is not hard to show that $v=f_i(H)$ in $C_i$, where
$f_i\in H^1_{loc}((0,1))$ for $i=1,3$ and $f_i\in H^1_{loc}((-1,0))$ for $i=2,4$.

It suffices to prove Lemma \ref{bd_l4} for $i=1$. The other cases are similar.

\medskip

\textbf{Step 1.} We first show that $f_1\in C^{2}((0,1))$.

In fact, let $\varphi(s)$ in $C^{\infty}_c(0,1)$. Since
\[
d\Delta v_A+|Dv_A|+AV(x)\cdot Dv_A=1
\]
Multiplying $\fee(H)$ on both sides and integrating by parts over the cell $C_1$ give
\[
\Rightarrow -d\int_{C_1}{Dv_A\fee'(H)DHdx}
+\int_{C_1}{|Dv_A|\fee(H)dx}
=\int_{C_1}{\fee(H)dx}
\]
Owing to Remark \ref{L1convergence}, sending $A\ra +\infty$,

\[
-d\int_{C_1}{Dv\fee'(H)DHdx}
+\int_{C_1}{|Dv|\fee(H)dx}
=\int_{C_1}{\fee(H)dx}
\]
Since $Dv=f'_{1}(H)DH$, by Coarea formula, it is easy to see that
$f_1=f_1(t)\in H^1_{loc}((0,1))$ is a weak solution of
\be\label{wellode}
d(f'_{1}(t)a(t))'+|f'_{1}(t)|b(t)=c(t)
\ee
where $a(t)$, $b(t)$, $c(t)>0$ and $\in C^{\infty}((0,1))$
$\Rightarrow f_1\in C^{2}((0,1))$.

\medskip

\textbf{Step 2.} We then prove that $f'_1< 0$. In fact, since
\[
\int_{C_1}{|v_A-v|^2dx}\ra 0\ \ \mbox{and}\ \
\int_{C_1}{|V(x)\cdot Dv_A|^2H^4dx}\ra 0,
\]
for any $0<a<b<1$, $\eps>0$, when $A$ is sufficiently large,
$\exists$ $0<a_{\eps}<b_{\eps}<1$ such that
$|a_{\eps}-a|<\eps$, $|b_{\eps}-b|<\eps$ and
\be\label{bd3}
\oint_{\left\{H=a_{\eps}\right\}}{|v_A-v|ds}\leq \eps\ \ ,\ \
\oint_{\left\{H=b_{\eps}\right\}}{|v_A-v|ds}\leq \eps
\ee
where $\oint$ denotes the integral average over a closed streamline, and
\be\label{bd4}
\max_{x,y\in\left\{H=a_{\eps}\right\}}|v_A(x)-v_A(y)|\leq \eps\ \ ,\ \
\max_{x,y\in\left\{H=b_{\eps}\right\}}|v_A(x)-v_A(y)|\leq \eps.
\ee
Since $d\Delta v_A+|Dv_A|+AV(x)\cdot Dv_A\geq 0$,
$v_A$ satisfies the maximum principle,
\[
 \max_{H=a_{\eps}}v_A\geq\max_{H=b_{\eps}}v_A
\]
According to (\ref{bd3}), (\ref{bd4}),
\[
f_1(a_{\eps})\geq f_1(b_{\eps})-4\eps
\]
Sending $\eps\ra 0$,
\[
 f_1(a)\geq f_1(b).
\]
Hence $f_{1}^{'}\leq 0$.  Owing to (\ref{wellode}), $f_{1}^{'}$ can not attain $0$.

\qed

The following lemma says that there is a mass loss of $|D v|$
as $A\to +\infty$. This implies the presence of
boundary layers (see Remark \ref{boundarylayer} and Figure \ref{internal}) where a positive
amount of mass of $|D v|$ is collected. More precisely,

\begin{lem}\label{loss}
$$
\lim_{\epsilon\to 0}\int_{\{|H|\geq \epsilon\}}|Dv|\,dx=\tau<1.
$$
\end{lem}

\nit {\bf Proof:} Owing to (\ref{L1bound}), it is obvious that $\tau\leq 1$.  Our goal is to exclude the case $\tau=1$.  We argue by contradiction. Let us assume that
\be\label{noloss}
\tau=1.
\ee
We first prove the following lemma.
\begin{lem}\label{bd_l5}
For $i=1,2,3,4$,
\[
\lim_{t\ra 0}f'_i(t)=0.
\]
\end{lem}

\nit {\bf Proof:} Since $d\Delta v_A+|Dv_A|+AV(x)\cdot Dv_A= 1$, for $a\in (0,1)$,
\[
\int_{\left\{|H|\geq a\right\}}{d\Delta v_Adx}
+\int_{\left\{|H|\geq a\right\}}{|Dv_A|dx}
=|\left\{|H|\geq a\right\}|.
\]
Here $|\mathcal{K}|$ represents the measure of the set $\mathcal {K}$. So
\[
d\int_{\left\{|H|=a\right\}}{\dpart{v_A}{n}ds}
+\int_{\left\{|H|\geq a\right\}}{|Dv_A|dx}
=|\left\{|H|\geq a\right\}|.
\]
Hence for $l\in (0,1)$ and small $\epsilon>0$,
\[
d\int_{\ell-\eps}^{\ell+\eps}da\int_{\left\{|H|= a\right\}}{\dpart{v_A}{n}ds}
+\int_{\ell-\eps}^{\ell+\eps}da\int_{\left\{|H\geq a\right\}}{|Dv_A|dx}
\]\[
=\int_{\ell-\eps}^{\ell+\eps}{|\left\{|H|\geq a\right\}|da}
\]
Sending $A\ra +\infty$, and by Remark \ref{L1convergence}, we deduce that
\[
\Rightarrow
d\int_{\ell-\eps}^{\ell+\eps}da\int_{\left\{|H|= a\right\}}{\dpart{v}{n}ds}
+\int_{\ell-\eps}^{\ell+\eps}da\int_{\left\{|H|\geq a\right\}}{|Dv|dx}
\]\[
=\int_{\ell-\eps}^{\ell+\eps}|{\left\{|H|\geq a\right\}|da}
\]
Dividing $\eps$ on both sides, we derive that
\be\label{green}
d\int_{\left\{|H|=a\right\}}{\dpart{v}{n}ds}
+\int_{\left\{|H|\geq a\right\}}{|Dv|dx}
=|\left\{|H|\geq a\right\}|
\ee
Since
\[
\lim_{a\ra 0}\int_{\left\{|H|\geq a\right\}}{|Dv|dx}
=\tau=1, \quad \text{by (\ref{noloss})}
 \]
\be\label{bd5}
\Rightarrow \lim_{a\ra 0}\int_{\left\{|H|=a\right\}}{\dpart{v}{n}ds}=0.
\ee
Note that
\[
\int_{\left\{|H|=a\right\}}{\dpart{v}{n}ds}
=|f'_1(a)|\int_{\left\{|H|=a\right\}\cap C_1}{|\nabla H|ds}
+|f'_2(-a)|\int_{\left\{|H|=-a\right\}\cap C_2}{|\nabla H|ds}
\]\be\label{bd6}
+|f'_3(a)|\int_{\left\{|H|=a\right\}\cap C_3}{|\nabla H|ds}
+|f'_4(-a)|\int_{\left\{|H|=-a\right\}\cap C_4}{|\nabla H|ds}.
\ee
Clearly for $i=1,2,3,4$
\[
\lim_{a\ra 0}\int_{\left\{|H|=a\right\}\cap C_i}{|\nabla H|ds}>0
\]
So Lemma \ref{bd_l5} holds.
\qed

Now let us finish the proof of Lemma \ref{loss}. By (\ref{green}), for any $0<a<b<1$,
\[
\Rightarrow
d\int_{\left\{|H|=a\right\}}{\dpart{v}{n}ds}
-d\int_{\left\{|H|=b\right\}}{\dpart{v}{n}ds}
+\int_{\left\{a\leq |H|\leq b\right\}}{|Dv|dx}
\]\[
=|\left\{a\leq |H|\leq b\right\}|
\]

According to (\ref{bd5}) and (\ref{bd6}),
\[
\lim_{a\ra 0}\int_{\left\{|H|=a\right\}}{\dpart{v}{n}ds}=0\ \ ,\ \
\int_{\left\{|H|=b\right\}}{\dpart{v}{n}ds}\geq 0
\]
\[
\Rightarrow
\int\!\!\!\!\int_{\left\{0<|H|\leq b\right\}}{|Dv|dx}
\geq |\left\{0< |H|\leq b\right\}|
\]
\[
\Rightarrow
1\leq{\int\!\!\!\int_{\left\{0<|H|\leq b\right\}}{|Dv|dx}
\over |\left\{0< |H|\leq b\right\}|}
=\sum_{i=1}^{4}{\int\!\!\!\int_{\left\{0<|H|\leq b\right\}\cap C_i}{|\nabla H||f_{i}^{'}(H)|dx}
\over |\left\{0< |H|\leq b\right\}|}
\longrightarrow 0
\]
as $b\ra 0$ according to Lemma \ref{bd_l5}. This is a contradiction.
\qed

\begin{rmk}\label{boundarylayer}(Presence of boundary layer)
Combining Remark \ref{L1convergence} and Lemma \ref{loss}, we have that
for any fixed $\epsilon>0$,
$$
\lim_{A\to +\infty}\int_{\{|H|\leq \epsilon\}}|Dv_A|\,dx= 1-\tau>0.
$$
\end{rmk}

Let $T$ be the smooth
solution of the following steady linear advection-diffusion problem
\be
d \Delta T+ AV(x)\cdot DT=0, \label{Teq}
\ee
subject to $T-e_1\cdot x$ being periodic and  $\int_{\Bbb T^2}T\,dx=0$. 
The following lemma says that the analysis of $\lam_A$ boils down to that of $L^1$ norm of $DT$.
\begin{lem}\label{tolinear}
There exists a constant $C$ such that
$$
{1\over C}||DT||_{L^1(\Bbb T^2)}\leq \lambda_{A} \leq C||DT||_{L^1(\Bbb T^2)}.
$$
In particular, 
\be\label{liminf}
\limsup_{A\to +\infty}{||DT||_{L^1(\Bbb T^2)}\over \lambda_A}\leq 1.
\ee
\end{lem}

\nit {\bf Proof}: We may assume that $v_A\rightharpoonup v$  in $L^{2}(\Bbb T^2)$ as $A \to +\infty$.  
Otherwise, we can argue by contradiction and use a subsequence.  
Let $S=G-T$ and $\beta_{A}\equiv \|DT\|_{L^1(\Bbb T^2)}$. Jensen's inequality and 
that $T - e_1\cdot x$ being periodic imply that $\beta_A \geq \|e_1\|_{L^1(\Bbb T^2)}=1$. 
The $S$ function satisfies the equation:
\be
d\Delta S+|DG|+AV(x)\cdot DS=\lam_A. \label{Seq}
\ee
Multiplying $S$ on both sides of (\ref{Seq}),
using integration by parts and the incompressibility of $V$,
we derive that
$$
d\int_{\Bbb T^2}|DS|^2\,dx=\int_{\Bbb T^2}S|DG|\,dx.
$$
A modification of Proposition 4 in \cite{NR} says that for a constant $C$
\[
G\leq C (1+||DG||_{L^1(\Bbb T^2)})=C(1+\lam_A)
\]
and 
$$
|T|\leq C\,(1+||DT||_{L^1(\Bbb T^2)})= C(1+\beta_{A}).
$$
Hence
$$
\int_{\Bbb T^2}|DS|^2\,dx\leq C (\lam_{A}^{2}+\beta_{A}^{2}).
$$
For $\epsilon>0$ which will be chosen later,
\be
\int_{\Bbb T^2}|DS|^2\,dx\geq \int_{\{|H|\leq \epsilon\}}|DS|^2\,dx\geq {1\over |\{|H|\leq \epsilon\}|}(\int_{\{|H|\leq \epsilon\}}|DS|\,dx)^2.
\label{fin1}
\ee
Note that
\be
\int_{\{|H|\leq \epsilon\}}|DG|\,dx \, \leq \,
\int_{\{|H|\leq \epsilon\}}|DS|\,dx + \int_{\{|H|\leq \epsilon\}}|DT|\, dx.
\label{fin2}
\ee
According to Remark \ref{boundarylayer}, when $A$ is large enough,
\be
\int_{\{|H|\leq \epsilon\}}|DG|\,dx\geq {(1-\tau)\over 2}\lambda_A. \label{fin3}
\ee
It follows from (\ref{fin1})-(\ref{fin3}) that:
\be
{(1-\tau)\over 2}\lambda_A \leq  C\, (\lambda_A^2 + \beta_{A}^{2})^{1/2} |\{|H|\leq \epsilon\}|^{1/2}
+ \beta_{A}. \label{fin4}
\ee
Since $\lim_{\epsilon\to 0}|\{|H|\leq \epsilon\}|=0$, we may choose $\epsilon$ small enough such that
the first term on the right hand side of (\ref{fin4}) is bounded from above by
$$
{(1-\tau)\over 4} (\lambda_A + \beta_A),
$$
implying:
\be
\lambda_A \leq C\beta_A. \label{fin5}
\ee
To finish the proof, it suffices to verify (\ref{liminf}). In fact, 
\ba
\lam_A & \geq & \int_{\{ |H|\leq \eps\}} |D G|\, dx  \no \\
& \geq & \int_{\{ |H|\leq \eps\}} |D T|\, dx -\left |\int_{\{ |H|\leq \eps\}} (|DG| - |DT|) \, dx \right | \no \\
& \geq & -\int_{\{ |H|\leq \eps\}} |DS |\, dx+ \int_{\{ |H|\leq \eps\}} |DT|\, dx \no \\
{\rm by} \; (\ref{fin1}) & \geq & - C(\lam_{A}^{2}+\beta_{A}^{2})^{1/2}\,|\{ |H|\leq \eps\}|^{1/2} +
 \int_{\{ |H|\leq \eps\}} |DT|\, dx. \no
\ea
For fixed $\delta>0$, we may choose $\epsilon$ sufficiently small such that
$$
C(\lam_{A}^{2}+\beta_{A}^{2})^{1/2}\,|\{ |H|\leq \eps\}|^{1/2}\leq \delta (\lambda_A+\beta_A).
$$
Then
\be
\lam_A (1+\delta) \geq - \delta \beta_A + \int_{\{ |H|\leq \eps\}} |DT|\, dx.\label{fin6}
\ee
According to (\ref{flat}) in the Appendix,  for fixed $(\epsilon,\delta)>0$, we can choose $A$ large enough such that
$$
\int_{\{|H|\geq \eps\}} |DT| \, dx \leq \delta \leq \delta \beta_A,
$$
where the last inequality is due to $\beta_A \geq 1$. It follows that
\[
\int_{\{ |H|\leq \eps\}} |DT|\, dx\geq (1-{\delta})\beta_A, 
\]
which implies that
$$
{\beta_{A}\over \lambda_A}\leq {1+\delta\over 1-2\delta}.
$$
Therefore
$$
\limsup_{A\to +\infty}{\beta_{A}\over \lambda_A}\leq {1+\delta\over 1-2\delta}.
$$
Then (\ref{liminf}) follows by sending $\delta\to 0$.
\qed

\nit {\bf {Proof of Theorem 2.1}}: Combining Lemma \ref{tolinear} and Lemma \ref{al} ($p=1$),  
we obtain the square root logarithm upper bound of $\lambda_A$. \qed

We shall present numerical evidence of the square root logarithm lower bound in section 4.

\section{Linear Law in Shear Flows}
\setcounter{equation}{0}
In this section,  we will investigate the front speed asymptotics for the shear
flow, i.e,  $V(x,y)=(v(y),0)$ where $v(y)$ is a smooth periodic function with mean zero, but not
identically zero.
%We also assume that $v$ is not constantly zero.
Unlike the cellular flow,  the turbulent flame speed from the shear flow
grows linearly with respect to A.  For the inviscid G-equation,
an explicit formula of $\bar H$ is given in \cite{EMS}.
Here we focus on the viscous G-equation.  We will also discuss the curvature dependent G-equation.

For $P=(m,n)$, the corresponding cell problem is reduced to an ODE
\be\label{1d}
-d\psi''+\sqrt {m^2+(n+\psi')^2}+A\, m\, v(y)=\lambda(A).
\ee
To simplify the notation, we write $\bar H(P,A,d)$ as $\lambda(A)$. If $m=0$,
it is obvious that $\lambda (A)=|n|$.  So throughout this section, we assume that $m\ne 0$.
We first show that the turbulent flame speed $\lambda (A)$ is enhanced
as $A$ increases.
\begin{theo}{} $\lambda= \lambda(A)$ as a function of $A\geq 0$ is convex and
strictly increasing.
\end{theo}
Proof:   The convexity follows immediately from the inf-max formula
$$
\lambda(A)=\inf_{\phi\in C^2(\Bbb T^1)}\max_{y\in  \Bbb T^1}\{-d\phi''+\sqrt {m^2+(n+\phi')^2}+Amv(y)\}.
$$
To prove that it is strictly increasing, it suffices to show that
$$
\lambda (0)=|P|<\lambda (A),  \quad \text{for all $A>0$}.
$$
This follows immediately from Jensen's inequality and the  strict convexity of the function  $f(t)=\sqrt {m^2+t^2}$.
\qed

Now let us look the asymptotic behavior of $\lambda(A) \over A$ as $A\to +\infty$.
Choose a solution $\psi$ (viscosity solution if $d=0$) of the cell problem with mean zero.  Denote $\psi_A={\psi\over A}$. Then $\psi_A$ satisfies that
$$
-d\psi_A''+\sqrt {{m^2\over A^2}+({n\over A}+\psi_{A}^{'})^2}+mv(y)={\lambda(A)\over A}.
$$
Since $\lambda (A)\over A$ is bounded, maximal principle implies that $\psi_{A}^{'}$ is uniformly
bounded.  Hence $\psi_A$ is equally continuous for both the inviscid ($d=0$) and
the viscous case ($d>0$). Upon a subsequence if necessary, we may assume
 that $(\psi_A, {\lambda(A)\over A})$ converges to $(\phi,\bar \lambda)$.
Stability of viscosity solutions implies that $(\phi,\bar \lambda)$ satisfies the following
cell problem:
\be\label{limit}
-d\phi''+|\phi'|+mv(y)=\bar \lambda
\ee
which is a special case of (\ref{1d}) for $P=(0,0)$ and $A=1$ subject to $\int_{\Bbb T^1}\phi\,dx=0$.
 Here $\Bbb T^1=[0,1]$. Therefore $\bar \lambda$ and $\phi'$ are uniquely given.  In particular, $\bar \lambda$ is positive.
Hence $\lambda (A)$ grows linearly for the shear flow.
\begin{theo}
$$
\lim_{A\to +\infty}{\lambda(A)\over A}=\bar \lambda=\bar \lambda (d)> 0.
$$
When $d=0$, $\bar \lambda(0)=\max_{\Bbb T^1}mv$.
\end{theo}
Proof We only need to show that $\bar \lambda>0$.  Taking integration on both
sides of (\ref{limit}) leads to
$$
\int_{0}^{1}|\phi'|\,dy=\bar \lambda.
$$
Since $v$ is not a constant, $\phi'$ can not vanish everywhere.  So $\bar \lambda$ must be positive.
\qed

Next we shall see how $\bar \lambda$ depends on the diffusivity constant $d$.
The following result says that the diffusion will slow down the
front propagation.
\begin{theo}  For $d>0$,  $\bar \lambda=\bar \lambda(d)$ is strictly decreasing as a function of $d$.\end{theo}
Proof: Let $w(y,d)=\phi'(y,d)$ and take derivative with respect to $d$ on both sides of (\ref{limit}). We get that
$$
-w'-dw_{d}^{'}+\mathrm{sign}(w)w_d=\bar \lambda_{d}.
$$
Let $h=-(dw_d+w)$.  We have that
\be\label{ode1}
h'-d^{-1}\mathrm{sign}(w)h=\bar \lambda_{d}+d^{-1}\,|w|.
\ee
If $\bar \lambda_{d}\geq 0$,
\be\label{ode2}
h'-d^{-1}\mathrm{sign}(w)h\geq 0.
\ee
Since $\int_{0}^{1}h\,dy=0$,  there exists $y_0$ such that $h(y_0)=0$.
According to (\ref{ode2}) and the periodicity of $h$,  we must have that
$$
h\equiv 0.
$$
Due to (\ref{ode1}),  $w\equiv 0$. So (\ref{limit}) implies that $v$ is a
constant function. This is a contradiction.\qed

Apparently,  $\lim_{d\to 0^{+}}\bar \lambda(d)=\bar \lambda (0)=\max_{\Bbb T^1}mv$.
A subsequent question is the convergence rate.
It is not obvious at all whether $\bar \lambda (d)$ is
differentiable at $d=0$ since the inviscid equation ($d=0$) has multiple
solutions and those solutions are not $C^1$.
In the following, we show that $\bar \lambda_d(0)=0$.
\begin{theo}\label{decay}
$$
\lim_{d \ra 0^{+}} \, {\bar \lambda (0)-\bar \lambda(d)\over d}=0.
$$
\end{theo}
Proof:  Without loss of generality,  we assume that $m=1$ and $v(0)=\max_{\Bbb T^1}v=0$. By Theorem 3.3,  $\bar \lambda (d)< \bar \lambda (0)=0$ for all $d>0$.

\nit Case 1:  Suppose that $0$ is the unique maximum point of $v$ in $[0,1)$. For $d>0$,  let $\phi=\phi(y,d)$ be the unique solution of
$$
-d\phi''+|\phi'|+v(y)=\bar \lambda
$$
satisfying  $\phi(0,d)=0$. Then
$$
\lim_{d\to 0}\phi=s(y)   \quad \text{uniformly in $\Bbb T^1$},
$$
where $s(y)$ is the unique viscosity solution of
$$
|s'|+v(y)=\bar \lambda(0)=0
$$
satisfying that $s(0)=0$ which is given by the formula
$$
s(y)=
\begin{cases}
\int_{0}^{y}(-v(t))\,dt  \quad \text{for $0\leq y\leq \bar y$}\\
\int_{y}^{1}(-v(t))\,dt \quad \text{for $\bar y\leq y\leq 1$}.
\end{cases}
$$
Here $\bar y\in (0,1)$  is the unique point which satisfies that
$$
\int_{\bar y}^{1}(-v(t))\,dt=\int_{0}^{\bar y}(-v(t))\,dt.
$$
Let
$$
\hat s(y)=
\begin{cases}
2\int_{0}^{y}(-v(t))\,dt \quad \text{for $y\geq0$}\\
2\int_{y}^{0}(-v(t))\,dt \quad \text{for $y\leq 0$}.
\end{cases}
$$
Choose $y_d\in [0,1]$ such that
$$
\phi(y_d,d)-\hat s(y_d)=\max_{\Bbb T^1}(\phi-\hat s).
$$
Since
$$
s(0)-\hat s(0)=0>s(y)-\hat s(y)  \quad \text{for $y\ne 0$},
$$
we have that
\be\label{zero}
\lim_{d\to 0}y_d=0.
\ee
Maximal principle implies that
$$
\phi''(y_d)\leq \hat s^{''}(y_d)\leq C|y_d|
$$
and
$$
\phi'(y_d)=\hat s^{'}(y_d).
$$
Hence
$$
0\geq \bar \lambda(d)=-d\phi''(y_d)+|\phi'(y_d)|+v(y_d)\geq -Cd|y_d|+|v(y_d)|\geq -Cd|y_d|.
$$
Accordingly, (\ref{zero}) implies Theorem \ref{decay}.

\nit Case II: $v$ has more than one maximum point.  Choose a smooth periodic function $L$ such that $L(0)=0$ and
$$
L(y)<0  \quad \text{for $y\in (0,1)$}.
$$
For $\delta>0$, write
$$
v_{\delta}=v+\delta L.
$$
Then $v_{\delta}$ has a unique maximum point at $y=0$ in $[0,1)$.  Let $\bar \lambda_{\delta}(d)$ be the corresponding asymptotic limit. Owing to the inf-max formula,
$$
\bar \lambda (d)\geq \bar \lambda_{\delta}(d).
$$
By case I,
$$
\lim_{d\to 0}{\bar \lambda_{\delta}(d)\over d}=0.
$$
Therefore Theorem \ref{decay} holds. \qed

\begin{rmk} It remains an interesting problem to study whether
$\bar \lambda (d)=\bar \lambda(0) + O(d^{r})$ for some power $r>1$. Computation
suggests that $r=2$. Table \ref{barLm} lists the raw and
$d^2$-scaled $\bar{\lambda}$ values
for $d \sim 0$, and suggests the quadratic behavior of $\bar{\lambda}$ in $d$.
The numerical values in table 1 are obtained from
the cell problem (3.2) with $m=1$, $v(y)=\cos(2\pi y)-1$, and so $\lambda (0)=0$.
We consider the time dependent problem
$$\phi_t-d\phi''+|\phi'|+(\cos(2\pi y)-1)=0\ ,\ \phi(t=0)=0.$$
Then $\bar \lambda $ is extracted from $\phi_t\ra -\bar \lambda$ uniformly
in $(0,1)$ as $t\ra +\infty$.

Spatial derivatives $\phi'$ and $\phi''$ are discretized by central differencing
with small enough grid size to ensure accuracy.
By symmetry, we have $\phi'(y)>0$ in $(0,1/2)$ and $\phi'(y)<0$ in
$(1/2,1)$ for any $t\geq 0$. Implicit Euler scheme can be readily used to
relax the time step constraint and speed up convergence to steady state.
The overall scheme is implicit in time, and second order in space.
\end{rmk}

The proof of Theorem \ref{decay} can be easily modified to recover a known interesting result in \cite{JKM}.  Denote
$$
\Gamma=\{y\in \Bbb T^1|\  v(y)=\max v\}=\{y_i\}_{i=1}^{m}.
$$
We assume that $v''(y_i)\ne 0$ for all $i$ and $|v''(y_i)|$ is strictly increasing  as $i$ varies from 1 to $m$.  Suppose that $\phi=\phi (y,d)$ is the unique solution of
$$
-d\phi''+{1\over 2}|\phi'|^2+v(y)=I_{d}
$$
subject to
$$
\int_{0}^{1}\phi\,dy=0.
$$
Then
\begin{theo}[Jauslin-Kreiss-Moser \cite{JKM}] Assume $\max_{\Bbb T^1}v=0$. Then
$$
\lim_{d\to 0} {I_d\over d}=-\sqrt {|v''(y_1)|}
$$
and $\phi$ uniformly converges to $\phi_0$ which is the periodic viscosity solution of
$$
{1\over 2}|\phi_{0}^{'}|^2+v(y)=0
$$
with a unique transition point $y_1$.
\end{theo}

\begin{rmk} Note that $\phi'$ and $\phi_{0}^{'}$ are solutions of the
viscous and inviscid Burgers equations respectively. The inviscid equation ${1\over 2}|\phi'|^2+v(y)=0$ has many solutions
even up to a constant when $v$ has multiple maximum points ($m>1$). A solution is uniquely determined by its
transition points, i.e, where it changes from decreasing to increasing.
The above theorem says that under some nondegeneracy conditions,
the vanishing viscosity method will select a unique  ``physical" solution.
Compared to the method in \cite{JKM}, our approach is more elementary and
can be easily extended to higher dimensions and more general Hamiltonians
at least when the Aubry set only consists of finitely many points
(see others approaches in \cite{AIPS} using stochastic control and random perturbation theories).
\end{rmk}
\begin{table}
\caption{Values of $\bar \lambda$ and $d^2$-scaled $\bar \lambda$ at small $d$ for $v(y)=\cos(2\pi y)-1$.} \label{barLm}
\begin{center}
\begin{tabular}{|c|c||c|}
\hline
$d$ & $-\overline{\lambda}$ & $-\overline{\lambda}/d^2$\\
\hline
4e-2 & 5.9414e-2 & 3.7134e+1 \\
\hline
2e-2 & 1.5540e-2 & 3.8850e+1\\
\hline
1e-2 & 3.9280e-3 & 3.9280e+1\\
\hline
4e-3 & 6.3428e-4 & 3.9642e+1\\
\hline
2e-3 & 1.5514e-4 & 3.8785e+1\\
\hline
1e-3 & 3.7061e-5 & 3.7061e+1\\
\hline

\end{tabular}
%\caption{Values of $\bar \lambda$ and $d^2$ scaled $\bar \lambda$ at small $d$ for $v(y)=cos(2\pi y)-1$.}
\end{center}
\end{table}

In combustion modeling, the laminar flame speed $s_{l}$ might also depend on the
curvature of the flame front. Peters \cite{Pet00} proposed the
following curvature dependent G-equation:
$$
G_t-d|DG|\mathrm{div}({DG\over |DG|})+s_l|DG|+V(x)\cdot DG=0.
$$
Here $\kappa=\mathrm{div}({DG\over |DG|})$ is the mean curvature of the flame front.
In general, we do not know how to prove the existence of the turbulent flame speed
for the curvature dependent G-equation. However, for the shear flow,
the corresponding cell problem is reduced to an ODE
$$
{-dm^2\phi''\over m^2+(n+\phi')^2}+\sqrt {m^2+(n+\phi')^2}+mv(y)=\lambda.
$$
Here we set $s_l=1$. It is very easy to verify the existence of classical solution
$\phi$ (unique up to an additive constant) and a constant $\lambda$.
Intuitively, the $\lambda$ from the curvature G-equation should be between the inviscid and the viscous case. The following theorem says that its asymptotic limit coincides with the inviscid case.

\begin{theo}\label{curvature}
If $v$ is scaled to $A\, v$, then $\lambda =\lambda(A)$ satisfies the growth law:
$$
\lim_{A\to +\infty}{\lambda (A)\over A}=\max_{\Bbb T^1}mv.
$$
\end{theo}
Proof:  Up to a subsequence if necessary, we may assume that
$$
\lim_{A\to +\infty}{\lambda (A)\over A}=\bar \lambda.
$$
Suppose that $\phi (y_0)=\min_{\Bbb T^1}\phi$. Then
$$
\lambda(A)\leq \sqrt {m^2+n^2}+Am v(y_0)\leq  \sqrt {m^2+n^2}+A\max_{\Bbb T^1}m v(y).
$$
Hence
$$
\bar \lambda\leq \max_{\Bbb T^1}mv.
$$
One the other hand, let $mv(y_1)=\max_{\Bbb T^1}mv$. Then for any $\delta>0$,
$$
\begin{array}{ll}
{\delta \lambda (A) \over A}&={1\over A}\int_{y_1}^{y_1+\delta}{-dm^2\phi''\over m^2+(n+\phi')^2}\,dy+{1\over A}\int_{y_1}^{y_1+\delta}(\sqrt {m^2+(n+\phi')^2}+Amv(y))\,dy\\[5mm]
&\geq -{dm\pi \over A}+\delta \min_{[y_1,y_1+\delta]}mv.
\end{array}
$$
Therefore $\bar \lambda\geq \max_{\Bbb T^1}mv$.  So Theorem \ref{curvature} holds.
\qed

\section{$\bar{H}$ in Cellular Flow: Empirical Law}
\begin{figure}\centering \label{barHA}
	\centering
\centerline{\includegraphics[width=4.5in,height=2.5in]{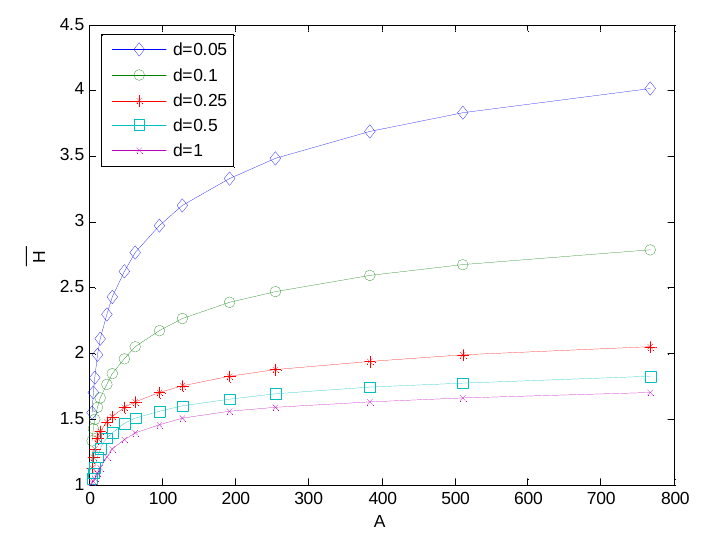}}
\caption{Plot of $\bar H(P,d)$ vs. $A \in [0,768]$ for $d=0.05,0.1,0.25,0.5,1$.}
\end{figure}

\setcounter{equation}{0}
Computation is carried out with finite difference discretization and iteration method
on equation (\ref{Ge4}) for $d \geq 0.1$
and upwind method on the evolution equation (\ref{Ge3}) for $d < 0.1$ with small enough grid size.
We choose $s_l=1$, $P=e_1$, and
$V(x)=(A/ 2\pi)\nabla^{\perp} \sin 2\pi x_1 \sin 2\pi x_2$. More details
can be found in \cite{LXY}. Numerical values of $\bar{H}(A,d)$ are
obtained for $d\in (0,1]$ and $A$ up to 768.
Figure \ref{barHA} clearly shows that $\bar H$ is decreasing with respect
to the diffusivity $d$ and increasing with respect to the flow intensity $A$. Such qualitative
property of $\bar{H}$ remains to be proved. The scaled quantity $\bar{H}/\sqrt{\log A}$ are
listed in Tables \ref{HlogA1}-\ref{HlogA2} for a range of decreasing $d$ values from 1 to 0.05.

We observe from Table \ref{HlogA1} that when $d\geq 0.5$,
the ratio ${\bar{H}\over\sqrt{\log A}}$ is slightly decreasing and approaching a limit
as $A$ increases. This is consistent with Theorem \ref{bd_t1}. For smaller $d$ values ($d\leq 0.1$), Table \ref{HlogA2} shows
that the ratio becomes slightly increasing as $A$ increases. Tables \ref{HlogA1}-\ref{HlogA2}
suggest the asymptotic law:
\be \label{conj}
\bar{H}(d,A) \sim c(d) \sqrt{\log A},\;\; d> 0 \;{\rm fixed}, \; A\gg 1,
\ee
 where $c(d)$ is a decreasing function of $d$.

It remains an interesting yet challenging problem to
analyze the empirical growth law (\ref{conj}), especially in the regime of small $d$.
\begin{table}
\caption{The $\sqrt{\log A}$ scaled effective Hamiltonian $\bar{H}(d,A)$ at large A and $d=1, 0.5, 0.25$.}\label{HlogA1}
\begin{center}
\begin{tabular}{|c||c|c||c|c||c|c|}
\hline
$A$& $\bar{H}_{(d=1)}$& ${\overline{H}\over\sqrt{\log A}}$& $\bar{H}_{(d=0.5)}$& ${\overline{H}\over \sqrt{\log A}}$&
$\bar{H}_{(d=0.25)}$& ${\overline{H}\over \sqrt{\log A}}$\\
\hline
  32&  1.2701&  6.8224e-1&  1.3987&  7.5132e-1& 1.5211&  8.1707e-1\\
  48&  1.3473&  6.8476e-1&  1.4634&  7.4377e-1& 1.5846&  8.0537e-1\\
  64&  1.3968&  6.8493e-1&  1.5058&  7.3838e-1& 1.6312&  7.9987e-1\\
  96&  1.4605&  6.8362e-1&  1.5606&  7.3047e-1& 1.7002&  7.9581e-1\\
 128&  1.5017&  6.8174e-1&  1.5975&  7.2524e-1& 1.7505&  7.9469e-1\\
 192&  1.5543&  6.7787e-1&  1.6496&  7.1943e-1& 1.8222&  7.9471e-1\\
 256&  1.5881&  6.7440e-1&  1.6869&  7.1636e-1& 1.8714&  7.9471e-1\\
 384&  1.6328&  6.6935e-1&  1.7393&  7.1300e-1& 1.9402&  7.9536e-1\\
 512&  1.6631&  6.6586e-1&  1.7741&  7.1030e-1& 1.9873&  7.9566e-1\\
 768&  1.7049&  6.6144e-1&  1.8227&  7.0714e-1& 2.0504&  7.9548e-1\\
\hline
\end{tabular}
\end{center}
\end{table}
\begin{table}
 \caption{The $\sqrt{\log A}$ scaled effective Hamiltonian $\bar{H}(d,A)$ at large A and $d=0.1, 0.05$.}\label{HlogA2}
\begin{center}
\begin{tabular}{|c||c|c||c|c|}
\hline
$A$& $\bar{H}_{(d=0.1)}$& ${\overline{H}\over \sqrt{\log A}}$& $\bar{H}_{(d=0.05)}$& ${\overline{H}\over \sqrt{\log A}}$ \\
\hline
  32&  1.8418&  9.8934e-1&  2.4290&  1.3048e+0\\
  48&  1.9592&  9.9576e-1&  2.6253&  1.3343e+0\\
  64&  2.0459&  1.0032e+0&  2.7698&  1.3582e+0\\
  96&  2.1696&  1.0155e+0&  2.9749&  1.3925e+0\\
 128&  2.2590&  1.0255e+0&  3.1242&  1.4183e+0\\
 192&  2.3839&  1.0397e+0&  3.3327&  1.4535e+0\\
 256&  2.4724&  1.0499e+0&  3.4795&  1.4776e+0\\
 384&  2.5938&  1.0633e+0&  3.6878&  1.5118e+0\\
 512&  2.6774&  1.0720e+0&  3.8292&  1.5331e+0\\
 768&  2.7900&  1.0824e+0&  4.0189&  1.5592e+0\\
\hline
\end{tabular}
\end{center}
\end{table}

\section{Conclusions}
We studied the front speed asymptotics in the viscous G-equation by analyzing
the related cell problem of homogenization. A new and striking
result is that for cellular flows and any positive viscosity in the viscous G-equation,
the front speed can not grow faster than $\sqrt {\mathrm{log} A}$. In contrast, the front speed
of the inviscid G-equation grows almost linearly
in the large amplitude regime of the cellular flows. In shear flows,
the front speed of the G-equation grows linearly in the large flow amplitude.
The growth rate is a monotone decreasing function of the viscosity coefficient.
The linear growth law in shear flows also
persists in the curvature dependent G-equation, with the same growth rate as that of
the inviscid G-equation.

\section{Appendix: $L^p$ bound of $|DT|$ for $p\in [1,2]$}
\setcounter{equation}{0}

Let $T$ be the smooth solution of the following steady diffusion-advection problem
$$
d \Delta T+ AV(x)\cdot DT=0
$$
subject to $T-e_1\cdot x$ being periodic and  $\int_{\Bbb T^2}T\,dx=0$.  It is known that
$$
\|DT\|_{L^2(\Bbb T^2)}\leq C(1+A^{1\over 4}).
$$
Also according to Theorem 4.2 in \cite{NPR},  $|DT|$ decays very fast away from those stream lines $\{H=0\}$. Precisely speaking, for $N \geq 1$:
\be\label{flat}
\|D T\|^{2}_{L^2(\{x \in \mathbb{T}: |H| \geq N \sqrt{\eps}\})} \leq {C\over \sqrt{\epsilon} N^4}. 
\ee
Here $A={1\over \epsilon}$.  Then we have that
\begin{lem}\label{al}
For $p\in [1,2]$,
$$
\|DT\|_{L^{p}(\Bbb T^2)}\leq CA^{{p-1\over 2p}}(\mathrm{log}A)^{{2-p\over 2p}}.
$$
Here $C$ is a constant independent of both $p$ and $A$.
\end{lem}

\nit The proof is a modification of the proof of Lemma 4.1 in \cite{NXY}.   Denote  $\Omega_{N}=\{x \in \Bbb T^2: (N-1) \sqrt{\epsilon}\leq |H| \leq N \sqrt{\epsilon}\}$. According to H\"older's inequality
$$
\begin{array}{ll}

\int_{\Bbb T^2}|DT|^{p}\,dx&=\sum_{N=1}^{1\over {\sqrt{\epsilon}}}\int_{\Omega_{N}}|DT|^{p}\,dx\\[5mm]
&\leq \sum_{N=1}^{1\over {\sqrt{\epsilon}}}{(\int_{\Omega_{N}}|DT|^{2}\,dx)}^{p\over 2}|\Omega_{N}|^{1-{p\over 2}}\\[5mm]
&\leq C {\epsilon}^{-p\over 4}|\Omega_{N}|^{1-{p\over 2}}(\sum_{N=1}^{\infty}(1+{1\over N^{s}})),

\end{array}
$$
where $s={2p}\geq 2$ and $|\Omega_N|$ is the measure of $\Omega_N$.  An easy computation shows that
$$
|\Omega_N|\leq C\sqrt {\epsilon} |\mathrm{log}(\epsilon)|.
$$
Hence the above lemma holds.  \qed

\begin{figure}
\centering
\centerline{\includegraphics[width=4.0in,height=2.5in]{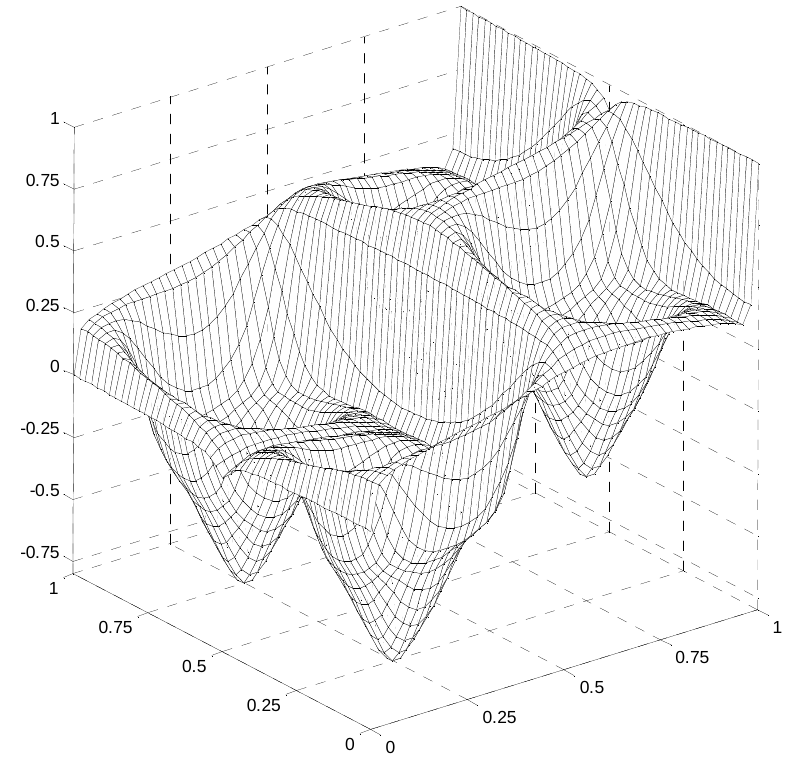}}
\caption{Graph of solution $G$ to equation (\ref{bd2}) in cellular flow at $A=16$, $d=0.01$:
presence of both internal wells and layers at quarter cell boundaries.} \label{internal}
\end{figure}
\begin{rmk} (Difference between $T$ and $G$)  As $A$ increases, it is known
that $T$ from equation (\ref{Teq}) becomes more and more like a
constant inside each quarter cell  \cite{NPR}.
This is not true for $G$ from (\ref{bd2}) due to the nonlinearity. In fact, according to Lemma \ref{bd_l4},  $v_A={G\over \lambda_A}$ converges to $v$ which is a strictly decreasing function of $|H|$ in each cell. Hence internal wells will emerge as ${A\over d} \gg 1$. This is confirmed by numerical calculations, see Figure \ref{internal}.

\end{rmk}

\bibliographystyle{plain}

\end{document}